\documentclass[12pt]{amsart}
\usepackage{amsmath,amsthm, amssymb}
\usepackage{amsrefs}
\usepackage[backref]{hyperref}

\theoremstyle{plain}
\newtheorem*{main theorem}{Main Theorem}
\newtheorem*{theorem*}{Theorem}
\newtheorem*{emptytheorem*}{}
\newtheorem{theorem}{Theorem}

\newtheorem*{proposition*}{Proposition}
\newtheorem{lemma}{Lemma}

\newtheorem{corollary}{Corollary}
\theoremstyle{definition}
\newtheorem{definition}{Definition}
\newtheorem*{definition*}{Definition}
\theoremstyle{remark}

\newtheorem*{remark*}{Remark}

\begin{document}
\title[Minimum principle for Lyapunov exponents]
{A minimum principle for Lyapunov exponents and a higher-dimensional 
version of a Theorem of Ma\~n\'e}

\author{Yongluo Cao}
\address{Department of Mathematics, 
Suzhou University, 
Suzhou 215006, Jiangsu, 
P.R. China}
\email{\href{mailto:ylcao@suda.edu.cn}{ylcao@suda.edu.cn}, 
\href{mailto:yongluo@pub.sz.jsinfo.net}{yongluo@pub.sz.jsinfo.net}}

\author{Stefano Luzzatto}
\address{Department of Mathematics, 
Imperial College, 
180 Queen's 
Gate,  London SW7 2AZ, UK}
\email{\href{Stefano.Luzzatto@imperial.ac.uk}
{Stefano.Luzzatto@imperial.ac.uk}}
\urladdr{\href{http://www.ma.ic.ac.uk/~luzzatto}
{http://www.ma.ic.ac.uk/\textasciitilde 
luzzatto}}

\author{Isabel Rios}
\address{Universidade Federal Fluminense, Niteroi, RJ, Brazil.}
\email{\href{mailto:rios@mat.uff.br}{rios@mat.uff.br}}

\thanks{IR was partially supported by CAPES and FAPERJ (Brazil). 
YC was partially supported by  NSF(10071055) and SFMSBRP of 
China and The Royal Society. The authors also wish to acknowledge the 
hospitality of Imperial College London where this work was carried out.}  
\date{22 August 2003}

\begin{abstract} 
We consider compact invariant sets 
\( \Lambda \) for \( C^{1}  \) maps in arbitrary dimension. 
We prove that if \( \Lambda \) contains no critical points then
there exists an invariant probability 
measure with a 
Lyapunov exponent \( \lambda \) which is the \emph{minimum} of 
all Lyapunov exponents for all invariant 
measures supported on \( \Lambda \).  We apply this result to prove 
that \( \Lambda \) is 
\emph{uniformly expanding}  if every invariant 
probability measure supported on \( \Lambda \) is hyperbolic 
repelling.  
This generalizes a well known theorem of Ma\~n\'e to the 
higher-dimensional setting. 
\end{abstract}

\maketitle

\section{Introduction and results}
Hyperbolicity as been a key idea in the modern theory of Dynamical 
Systems and a basic related problem has been that of proving hyperbolicity 
of particular systems or general classes under a priori weaker 
assumptions. In 1985, Ma\~n\'e \cite{Man85} 
proved a remarkable result in the 
setting of one-dimensional maps to the effect that any compact 
invariant set not containing any critical points and with all periodic 
points hyperbolic repelling is actually uniformly hyperbolic. The 
point here is that the periodic points are not assumed to be 
\emph{uniformly} hyperbolic, and that even if they were, this 
hyperbolicity does not necessarily extend to the whole invariant set. 
Generalizations of this result to the higher dimensional setting are 
quite problematic and there has been no substantial progress to date. 
Here we solve this problem under the slightly stronger but natural 
assumptions that all invariant measures are hyperbolic repelling, i.e. 
have positive Lyapunov exponents. As an intermediate result of 
independent interest we get that if all Lyapunov exponents are positive then 
they are uniformly positive and  the minimum of 
such exponents is actually realized for some invariant measure.

Throughout the paper, we let \( M \) be a compact Riemannian 
manifold of dimension \( d\geq 1 \) and let  \( f: M \to M \) 
be \( C^{1} \) map .  We say that \( x \) is 
a critical point for \( f \) if \( \det Df_{x} = 0 \). For a 
compact invariant set \( \Lambda \) we let 
\( \mathcal M(f)=\mathcal M(f, \Lambda) \) denote the set of all \( f 
\)-invariant Borel probability measures with support in \( \Lambda \) 
and \( \mathcal E(f, \Lambda) \) the subset of all ergodic invariant 
measures. 

\subsection{Lyapunov exponents}

\begin{definition}
    We say that \( \lambda \) is a \emph{Lyapunov exponent} for \( f \) 
    if there exists a point \( x \) and a vector \( v\in T_{x}M \) 
    such that 
    \[ 
    \lambda = \lim_{n\to\infty} \frac{1}{n}\log \|Df_{x}^{n}(v)\|.
    \]
    We let \( \mathcal L(f) \) denote the set of all Lyapunov 
    exponents for \( f \). 
\end{definition}

Classical results \cites{Ose68, Rue82} imply that 
for each ergodic measure \( \mu\in\mathcal M(f) \) 
 there exists a number \(  0< k\leq d \), constants \( 
    \lambda^{1}<\ldots<\lambda^{k} \), and a 
    \( \mu \)-measurable splitting \( T_{x}M= 
    E^{1}\oplus \ldots \oplus E^{k} \) of the tangent 
    bundle over \( \Lambda \), such that 
\(
    \lim_{n\to \infty} \frac{1}{n}\log\|Df^{n}_{x}(v)\|=\lambda^{j}
\)
    for \( \mu \)-almost every \( x \) and every non-zero 
    vector  \( 
    v\in E^{j}_{x}\subset T_{x}M \).  
 For non-ergodic measures 
 the number \( k \), the constants \( \lambda^{j} \), 
    and the tangent bundle decomposition may depend on the ergodic 
    component. 
The constants \( \lambda^{j} \) are called the Lyapunov 
exponents associated to the measure \( \mu \). 

\begin{definition}
We let \( \mathcal 
L(\mu) \) denote the set of all Lyapunov exponents associated to a 
given ergodic measure \( \mu \)
\end{definition}
Notice that exceptional points for all invariant measures may still 
have some have a well defined Lyapunov exponent which may be 
unrelated to any of the Lyapunov exponents of any ergodic invariant 
measure. Therefore in general we have 
 \[ 
    \bigcup_{\mu\in\mathcal E(f,\Lambda)} 
    \!\!\!\! \mathcal L(\mu) \subsetneqq \mathcal L(f).
 \]

\subsection{The minimum principle}
Lyapunov exponents do not in general vary continuous either with the 
base point \( x \) or with the measure, and so it is hard if not 
impossible to formulate any general compactness statements about the 
sets \( \mathcal L(f) \) of \( \cap_{\mu} \mathcal L(\mu) \).
Here we prove the following
\begin{theorem}[Minimum principle for Lyapunov exponents]
    \label{minimum principle}
 Suppose \( \Lambda \) does not contain any critical points. 
 Then there exists an ergodic measure \( \mu\in \mathcal M(f,\Lambda) 
 \) such that 
\[ 
\inf\{\mathcal L(f)\}\in \mathcal L(\mu).
\]
 \end{theorem}
This can be interpreted as a mild compactness result on
\( \mathcal L(f) \): 
there exists some invariant probability measure \( \mu \) with 
support in \( \Lambda \) with an associated Lyapunov exponent which 
realizes the infimum over all Lyapunov exponents associated  
to all invariant measures. It also says that no ``freak'' exponent can 
be less than a ``proper'' exponent associated to some ergodic 
invariant measure.
Exactly the same argument also shows that 
there exists some (other) measure \( \mu \) for which
\[ 
\sup\{\mathcal L(f)\}\in\mathcal L(\mu).
\]
We do not show however that \( \mathcal L(f) \) is closed. We 
conjecture this to be the case in this setting and believe that the 
absence of critical points is a necessary condition for such property 
to hold. 

In the application to be given below we 
shall consider the situation in which all invariant measures 
 have only positive Lyapunov exponents. The minimum principle 
then immediately implies the following 
\begin{corollary}
      Suppose \( \Lambda \)  contains no critical 
      points and 
      \[
      \mathcal L(f) > 0.
      \] Then there exists a constant \( 
 \lambda>0 \) such that 
 \[ 
 \mathcal L(f) > \lambda.
 \] 
 \end{corollary}
 
\subsection{Ma\~n\'e's Theorem in arbitrary dimension} 
Our next result says that in fact an even stronger statement holds in 
this case. before stating it we recall some standard notation.  
The set \( \Lambda \) is said to be 
\emph{uniformly expanding} if there exist constants \( C, 
\lambda > 0 \) such that 
\begin{equation}\label{uniform expansion} 
\| Df^{n}_{x} (v) \| \geq Ce^{\lambda n} \|v\|
\end{equation}
for all \( x\in \Lambda \), 
all \( v\in T_{x}M \) and all \( n\geq 1 \). 
A measure \( \mu \) is 
\emph{expanding} if \( \mathcal L(\mu) > 0 \). 
The definition of uniform hyperbolicity 
implies in particular that all invariant measures 
supported on \( \Lambda \) are ``uniformly'' expanding: there exists 
some \( \lambda > 0 \) such that \( \mathcal 
L(f) \geq \lambda > 0 \). The converse however is non-trivial. 

\begin{theorem}\label{uniform}
Suppose \( \Lambda \)  contains no critical points and 
\[
\mathcal L(f) > 0. 
\]
 Then  \( \Lambda \)   is uniformly expanding. 
 \end{theorem}  

This is a generalization to the higher-dimensional setting of the well 
known theorem of Ma\~n\'e \cite{Man85} which gives the same 
conclusions in the one-dimensional setting 
under the weaker assumption that every periodic point is expanding. 
This has become almost a Folklore Theorem in One-Dimensional dynamics 
for the fundamental 
role it plays in a huge number of arguments. By comparison, the 
theory of higher-dimensional non-uniformly expanding maps is just 
taking off, see \cites{Alv03, Luz03}, and, by analogy with the 
one-dimensional case, we feel that our results might play some 
significant role in that theory. 
 
We emphasize that the absence of critical points in \( \Lambda \) is 
a necessary condition. There are many quite generic
situations, including for example Collet-Eckmann 
one-dimensional maps
\cite{NowSan98}, of 
 compact invariant sets \( \Lambda \) 
satisfying \( \mathcal L(f) \geq \lambda > 0 \) but which 
 contain a critical point and thus  clearly cannot be uniformly 
 expanding. Partial results in the direction of our 
 Theorem \ref{uniform} were obtained in \cites{AlvAraSau, Cao03} in the 
 context of globally \( C^{1} \) local diffeomorphisms. Here we 
 provide a generalization of those results by using a quite 
 different argument based on the minimum principle given above.
Theorem \ref{uniform} also implies  the following  
 \begin{corollary}\label{uniform corollary}
     Let \( f: M \to M \) be a \( C^{1+\alpha} \) local diffeomorphism 
     of a compact Riemannian manifold \( M \) and suppose that 
     \[ 
     \mathcal L(f) > 0.
     \]
     Then there exists a (unique) measure \( \mu\in\mathcal E(f) \) 
     which is absolutely continuous with respect to the Riemannian 
     volume on \( M \). 
 \end{corollary}     
 The conclusion follows from the uniform hyperbolicity of \( f \) by 
 well known and absolutely classical arguments. However we want to 
 emphasize here that the assumptions differ from most kinds of 
 assumptions used to imply the existence of absolutely continuous 
 invariant measures, both in the uniform and non-uniform setting, in at 
 least one important feature. They do not, \emph{a priori}, say 
 anything about the dynamical or hyperbolic properties of 
 a positive measure set of points for the Riemannian volume. They just 
 specify that any invariant measure must have positive Lyapunov 
 exponent. A priori \emph{all} of these measures may be singular with respect 
 to the Riemannian volume. 
 Notice that for this corollary we need the derivative of \( f \) to 
 be H\"older continuous, since the proof of the existence of an 
 absolutely continuous invariant measure requires distortion 
 estimates which require at more regularity than that 
 needed for the other results.
 
 \subsection{Fibred maps}
 Some  straightforward adaptations of the definitions and 
 the arguments allow us to obtain our results in the 
 more general setting of fibred maps. More precisely, we assume as 
 above that \( f: M \to M \) is a \( C^{1} \) map of a compact Riemannian 
 manifold and that \( \Lambda \) is a compact invariant set. We now 
 assume that the tangent bundle \( T_{\Lambda}M \) over \( \Lambda \) 
 admits a continuous decomposition \( T_{\Lambda}M= 
 E^{1}_{\Lambda}\oplus E^{2}_{\Lambda} \) into \( Df \)-invariant 
 subbundles. In particular the angles between the subspaces \( 
 E^{1}_{x} \) and \( E^{2}_{x} \) are uniformly bounded below for all \( 
 x\in\Lambda \). The (measurable) Oseledets-Ruelle 
 decomposition must be consistent with this continuous decomposition and 
 therefore it makes sense to talk about the sets \( \mathcal 
 L^{1}(\mu) \) and \( \mathcal L^{2}(\mu) \) of Lyapunov exponents 
 in the directions of \( E^{1} \) and \( E^{2} \) respectively. 
 Notice that \( \mathcal L(\mu) = \mathcal L^{1}(\mu) \cup \mathcal 
 L^{2}(\mu) \). We then have the following
 
\begin{theorem}\label{fibre min prin}
 Suppose \( \Lambda \) does not contain any critical points. 
 Then
\[ 
\inf\{\mathcal L^{1}(f)\}\in \mathcal L^{1}(f).
\]
 \end{theorem}
 The definition of uniform hyperbolicity given above 
 can also be generalized and we say that the set \( \Lambda \) is 
\emph{uniformly expanding in the direction of } \( E^{1} \) 
if there exist constants \( C, 
\lambda > 0 \) such that 
\begin{equation}\label{fibre un ex} 
\| Df^{n}_{x} (v) \| \geq Ce^{\lambda n} \|v\|
\end{equation}
for all \( x\in \Lambda \), 
all \( v\in E^{1}_{x}\) and all \( n\geq 1 \).
 Then we have the following
 
 \begin{theorem}\label{fibre uniform}
Suppose \( \Lambda \)  contains no  critical points and 
 \[
 \mathcal L^{1}(f) > 0.
 \]
 Then  \( \Lambda \)   is uniformly expanding in the direction of \( E^{1} \).  
 \end{theorem}  

 We emphasize that we are not assuming here any \emph{partial 
 hyperbolicity}. The results do not depend in anyway on the 
 hyperbolicity properties in the complementary subspace to the one 
 under consideration. Theorems \ref{minimum principle} and 
 \ref{uniform} are special cases of Theorems \ref{fibre min 
 prin} and \ref{fibre uniform} corresponding to the case in 
 which the complementary subbundle \( E^{2}_{\Lambda} \) is trivial. 
The proofs of the more general cases differ from the special cases 
only in the notation and so, in order to keep this
 as simple as possible, we shall prove the results explicitly in the
 special cases.

 \subsection{Remarks}
Before starting the proof of the Theorems, we make some observations 
concerning the setup. In particular we want to emphasize the 
difference between the uniformity statement in Corollary 1 and that in 
Theorem \ref{uniform}. 

An  invariant measure \( \mu \) has at most \( d \) 
distinct Lyapunov exponents and thus the condition \( \mathcal L(\mu)>0 \) 
implies that there exists \( \lambda^{1}(\mu) = \inf\{\mathcal L(\mu)\} 
\). By the definition of Lyapunov exponents this implies 
    \begin{equation}\label{liminf} 
    \liminf_{n\to \infty} \frac{1}{n}\|Df^{n}_{x}(v)\|\geq \lambda^{1} 
    \end{equation}
for \( \mu \)-almost every \( x \) and for every non-zero \( v\in 
T_{x}M \) (the actual limit exists only for those vectors in the 
appropriate subspace \( E^{1}_{x}\subset T_{x}M \)).     
In particular (directly from the definition of \( \liminf 
    \)) there exists, for every \( \lambda^{1} > \lambda > 0 \), 
    a measurable function \( C_{x} \), non-zero \( \mu \)-almost 
    everywhere, such that 
    \begin{equation}\label{expansion} 
    \|Df^{n}_{x}(v)\|\geq C_{x}e^{\lambda n} \|v\|
    \end{equation}
    for every non-zero vector \( v\in T_{x}M \) and every \( n\geq 1 
    \). 
    Crucially here, the liminf in \eqref{liminf} cannot be assumed 
    to be achieved uniformly and thus the constant \( C_{x} \) is in 
    general not uniformly bounded away from \( 0 \).  
    
    Thus, the step from \( \mathcal L(f)>0 \) to uniform expansivity 
    requires two uniformity estimates to be obtained. The first, given 
    by Corollary 1, implies that 
    \eqref{expansion} holds on a set of total probability, 
    (i.e. on a set \( \mathcal B \) such that \( \mu(\mathcal B)=1 \) for all 
    measures \( \mu\in\mathcal M(f) \)) with a \emph{uniform bound on the 
        growth rate} \( \lambda \). 
	However, this still leaves us with a family of \emph{measurable} 
    functions \( C_{x}(\mu) \) which are not, a priori,
    uniformly bounded below. 
    Theorem \ref{uniform} says that in the 
    absence of critical points they are indeed uniformly bounded 
below. Therefore the inequality \eqref{expansion} can be 
    extended to every point of \( \Lambda \) for uniform constants \( 
    \lambda \) and \( C \), as required in the definition 
    \eqref{uniform expansion} of uniform expansion. 
    
After writing this paper we became aware of some previously published 
papers \cites{JohPalSel87, Slo95, Slo97,StaStu00} containing related 
results concerning uniformity of Lyapunov exponents and other 
quantities in the spirit of our Minimum Principle. The arguments given 
there are significantly more general and the proofs significantly more 
complicated and it is not clear that our result follows immediately 
from the statements in any of these papers without some additional 
non-trivial arguments. We thank Jaroslav Stark for bringing these 
references to our attention.

 \section{Setup and notation}
 
 A key tool in our approach is that of lifting certain quantities to 
 the unit tangent bundle 
 \[
 SM= S_{\Lambda}M =  \{(x,v) \in TM: x\in \Lambda, v\in T_{x}M, 
 \|v\| =1 \}.
 \]
 Notice that \( SM \) is a compact, metric, measure space. 
 We let \( \pi: SM \to M \) denote the standard projection \( \pi (x,v) = x \). 
 We start by defining the map \( F: SM \to SM \) by 
 \[ 
 F (x,v) = \left(f(x), \frac{Df_{x}(v)}{\|Df_{x}(v)\|}\right).
 \]
 Notice that \( F \) is well defined on \( \Lambda \) since 
 \( \Lambda \) contains no critical points and thus \( 
 \|Df_{x}(v)\| \neq 0 \). We define iterates of \( F \) by 
 \[ 
 F^{n}(x,v) = \left(f^{n}(x), 
 \frac{Df^n_{x}v}{||Df^n_{x}v||} \right).
 \]
 We let \( \mathcal M(F) \) denote the space 
 of \( F \)-invariant probability measures on \( SM \) and let 
 \[ 
 \pi^{*}: \mathcal M(F) \to \mathcal M(f) 
 \]
 denote the standard 
 projection of measures where \( \pi^{*}\mu (A) = \mu (\pi^{-1}(A)) \) 
 for any \( \mu\in \mathcal M(F) \) and Borel set \( A\subset \Lambda 
 \). Let \( \mathcal E(F) \subset \mathcal M(F) \) and \( \mathcal 
 E(f)\subset \mathcal M(f) \) denote the subsets of ergodic invariant 
 measures for \( F \) and \( f \) respectively. 
 
 \begin{lemma}\label{pistar}
     The projection \( \pi^{*} \) sends \( \mathcal E(F)\) to 
     \( \mathcal E(f) \)
     and the restriction \( \pi^{*}: \mathcal E(F) \to \mathcal E(f) \)
     is surjective.
 \end{lemma}
 
 \begin{proof} 
     We show first of all that \( \pi^{*}(\mathcal E(F))\subseteq 
     \mathcal E(f) \). Let \( \mu\in \mathcal E(F) \) and let \( 
     \mu^{*}=\pi^{*}\mu \). Suppose that \( A 
     \subset M \) satisfies \( f^{-1}(A)=A \). Then  \( 
     F^{-1}(\pi^{-1}(A))=\pi^{-1}(A) \) and by ergodicity we have 
     \( \mu(\pi^{-1}(A) ) =\) 0 or 1.
Then by the definition of \( \pi^{*} \) we have 
     \( \mu^{*}(A)=\pi^{*}\mu(A) = \mu (\pi^{-1}(A)) =\) 0 or 1.
     Thus \( \mu^{*} \) is ergodic. 
     
     We now show that \( \pi^{*}:\mathcal M(F) \to \mathcal E(f) \) is 
     surjective by fixing a measure \( \mu^{*}\in\mathcal E(f) \) and 
     finding a measure \( \mu\in\mathcal M(F) \) with \( \pi^{*}\mu = \mu^{*} \). 
 By the Birkhoff 
     ergodic theorem we know that the set 
     \[ 
     X = \left\{x: \frac{1}{n}\sum_{i=0}^{n-1}\delta_{f^{i}(x)} \to 
     \mu^{*}\right\} 
     \ \text{ satisfies } \ \mu^{*}(X)=1 
     \]
For such a generic point \( 
     x\in X \) and some vector \( v\in T_{x}M \), consider the sequence of 
     measures 
     \[ 
     \mu_{n}= \frac{1}{n}\sum_{i=0}^{n-1}\delta_{F^{i}(x, v)}
     \]
     and let \( \mu\in\mathcal M(F) \)  be the limit of some 
     subsequence  \( \mu_{n_{k}} \). 
     We claim that \( \pi^{*}\mu=\mu^{*} \). 
     To see this, consider a continuous test function \( \varphi: M 
     \to \mathbb R \). We claim that 
     \[ 
     \int_{M}\varphi d(\pi^{*}\mu) = \int_{SM}\varphi \circ \pi 
     d\mu = \int_{M}\varphi d\mu^{*}.
     \]
     Since \( \varphi \) is chosen arbitrarily this implies the 
     claim. The first equality follows immediately. To obtain the 
     second we write, for \( x\in X \) and \( v\in T_{x}M \) as 
     above,  
     \begin{align*}
     \int_{M}\varphi d\mu^{*}& = 
     \lim_{n\to\infty}\frac{1}{n}\sum_{i=0}^{n-1} \varphi (f^{i}(x)) 
     \\ 
     &=   \lim_{k\to\infty}\frac{1}{n_{k}}\sum_{i=0}^{n_{k}-1} \varphi (f^{i}(x)) 
     \\
  &= 
      \lim_{k\to\infty}\frac{1}{n_{k}}\sum_{i=0}^{n_{k}-1} 
      (\varphi\circ \pi) (F^{i}(x, v)) \\ 
      &= 
      \int_{SM} (\varphi\circ \pi) d\mu.
     \end{align*}
The measure  \( \mu \) is not necessarily ergodic, but we claim that any ergodic 
component \( \tilde\mu \) of \( \mu \) also satisfies \( 
\pi^{*}\tilde\mu=\pi^{*}\mu=\mu^{*} \). Indeed consider the set \( 
X\subset M \) of \( \mu^{*} \) generic points defined above. Since 
\( \mu^{*}(X)=1 \) we also have that \( \mu(\pi^{-1}(X))=1 \) and 
therefore also \( \pi^{*}\tilde\mu(X) = \tilde\mu (\pi^{-1}(X)) =1 
\). 
Moreover, \( \pi^{*}\tilde\mu \) is ergodic by the first part of the 
statement in the Lemma, and therefore the Dirac averages of almost 
every point in \( X \) converge to \( \pi^{*}\tilde\mu \) implying 
that \( \pi^{*}\tilde\mu = \mu^{*}. \)
 \end{proof}
    
  \section{Minimum principle}
In this section we prove Theorem \ref{minimum principle}. 
We  define the observable \( \phi: SM \to \mathbb R \) by 
 \[ 
 \phi (x, v) = \log \|Df_{x}(v)\|.
 \]
 Notice that \( \phi \) is \emph{continuous} and 
 that we have 
 \begin{lemma}
     For every \( (x,v)\in SM \) and every \( n\geq 1 \) we have 
     \[ 
     \frac{1}{n}\log  \|Df_{x}^{n}(v)\| = 
     \frac{1}{n}\sum_{i=0}^{n-1}\phi (F^{i}(x,v)).
     \]
     In particular 
     \[ 
     \lambda (x,v) = \lim_{n\to\infty} 
       \frac{1}{n}\log  \|Df_{x}^{n}(v)\| = \lim_{n\to\infty} 
     \frac{1}{n}\sum_{i=0}^{n-1}\phi (F^{i}(x,v))
     \]
     whenever such limits exist (and the existence of one limit 
     implies the existence of the other).
  \end{lemma}
  \begin{proof} 
      Write first of all 
      \[ 
      \begin{split}\frac{\|Df^{n}_{x}(v)\|}{\|v\|} = 
       \frac{\|Df_{f^{n-1}(x)}(Df^{n-1}_{x}(v))\|}{\|Df^{n-1}_{x}(v)\|} 
        \frac{\|Df_{f^{n-2}(x)}(Df^{n-2}_{x}(v))\|}{\|Df^{n-2}_{x}(v)\|}
	\ldots \\ \ldots
 \frac{\|Df_{f(x)}(Df_{x}(v))\|}{\|Df_{x}(v)\|} 
  \frac{\|Df_{x}(v)\|}{\|v\|}. 
  \end{split}
      \]
Then, taking logs and using the definition of \( \phi \) we have 
     \begin{align*}
     \log \frac{\|Df^{n}_{x}(v)\|}{\|v\|}  &= 
     \sum_{i=0}^{n-1}\log \frac{\|Df_{f^{i}(x)} 
     (Df^{i}_{x}(v))\|}{\|Df^{i}_{x}(v)\|} \\ 
     &= 
      \sum_{i=0}^{n-1}\log \left \| Df_{f^{i}(x)} 
     \left(\frac{Df^{i}_{x}(v)}{\|Df^{i}_{x}(v)\|}\right)\right\|
    \\ & = 
     \sum_{i=0}^{n-1} \phi (F^{i}(x,v)).
     \end{align*}
     \end{proof}

 \begin{lemma}\label{minimumlemma}
     There exists a measure \( \hat\mu \in \mathcal E(F) \)
     such that 
     \[ 
     \int_{SM}\phi   d\hat\mu = \inf_{\mu\in\mathcal 
     M(F)}\int_{SM} \phi  d\mu.
     \]
     \end{lemma}
     \begin{proof}
By the compactness of \( \mathcal M(F) \) and continuity of \( \phi 
\) and of the integral functional \(
\ell_{\phi} (\mu) = \int  \phi d\mu,
\)
it follows that there exists some measure \( \bar\mu \) for which the the 
equality in the statement of the Lemma holds. This measure is not 
necessarily ergodic but we claim that some ergodic component of \( 
\bar\mu \) also satisfies the required equality. Indeed, by the Ergodic 
Decomposition Theorem \cite{Wal82}*{page153}, there exists a measure \( \tau \) on \( 
\mathcal M(F) \)  and a set 
\( \mathcal E_{\bar\mu}(F)\subseteq\mathcal 
E(F) \) with \( \tau(\mathcal E_{\bar\mu}(F))=1 \)  such that we have 
\( 
\int_{\mathcal E(F)} \nu d\tau = \bar \mu 
\)
in the sense that
\begin{equation}\label{measure1} 
\int_{\mathcal E(F)} \left(\int_{SM}\phi d\nu\right) d\tau = 
\int_{SM}\phi d\bar\mu. 
\end{equation}
Now, for any \( \nu\in\mathcal E_{\bar\mu}(F) \) we have 
\begin{equation}\label{measure2} 
\int_{SM} \phi d\nu \geq \inf_{\mu\in\mathcal 
     M(F)}\int_{SM} \phi  d\mu = 
     \int \phi d\bar\mu.
\end{equation}
If the inequality was strict for a positive \( \tau \) measure set, 
we would have 
\[ 
\int_{\mathcal E(F)} \left(\int_{SM}\phi d\nu\right) d\tau 
> \int \phi d\bar\mu
\]
contradicting \eqref{measure1}. Therefore there must be a \( \tau 
\) full measure (in particular non-empty)
set  of  measures in \( \mathcal E_{\bar\mu}(F) \) for which 
the equality  in \eqref{measure2} holds. 
\end{proof}

\begin{lemma}\label{lemmalyap1}
    \( \forall \ \lambda\in\mathcal L(f) \) 
    \( \exists \ \mu\in\mathcal M(F) \) such that 
    \[
    \int \phi d\mu = \lambda.
    \]
    \end{lemma}
    \begin{proof}
By assumption, there exists some \( (x,v)\in SM \) such that 
\begin{equation}\label{lyap}
\lambda(x,v) = \lim_{n\to\infty} \frac{1}{n}\log {\|Df^{n}_{x}(v)\|} 
= \lambda.  
\end{equation}
For such a point \( (x,v) \), consider the sequence of measures
\[ 
\mu_{n}=\frac{1}{n}\sum_{i=0}^{n-1}\delta_{F^{i}(x,v)}
\]
and let \( \mu\in\mathcal M(F) \)  be the limit of some subsequence 
 \( \mu_{n_{k}}\).  For any \( n\geq 1 \) we have 
\[ 
\frac{1}{n}\sum_{i=0}^{n-1}\phi(F^{i}(x,v)) = \int \phi d\mu_{n}.
\]
Therefore, by the definition of weak-star convergence, 
\[ 
\lim_{k\to\infty}
\frac{1}{n_{k}}\sum_{i=0}^{n-1}\phi(F^{i}(x,v)) = 
\lim_{n_{k}\to\infty}\int \phi d\mu_{n} 
= \int \phi d\mu.
\]
Moreover, using the existence of 
the limit in \eqref{lyap}, we have
\[ 
\lim_{k\to\infty} \frac{1}{n_{k}}\log 
\frac{\|Df^{n}_{x}(v)\|}{\|v\|}=\lim_{n\to\infty} \frac{1}{n}\log 
\frac{\|Df^{n}_{x}(v)\|}{\|v\|} = \lambda (x,v) = \lambda. 
\]
Finally, the definition of \( \phi \) implies 
\[ 
\lim_{k\to\infty} \frac{1}{n_{k}}\log 
\frac{\|Df^{n}_{x}(v)\|}{\|v\|}=
\lim_{k\to\infty}
\frac{1}{n_{k}}\sum_{i=0}^{n-1}\phi(F^{i}(x,v)).
\]
Substituting this into the two previous equations completes the proof.
\end{proof}

\begin{lemma}\label{lemmalyap2}
 \( \forall \ \mu\in\mathcal E(F) \)  \( \exists \ \lambda\in\mathcal L(f) \) 
 such that 
 \[
 \lambda = \int\phi d\mu.
 \]
  \end{lemma}  
  \begin{proof}
      By the ergodicity of \( \mu \) there exists a set \( A\subseteq 
      SM \) with \( \mu(A) = 1 \) such that for all \( (x,v)\in A \) 
      there exists a constant 
      \begin{equation}\label{lyap2} 
      \lambda=\lambda (x,v) = \lim_{n\to\infty} \frac{\|Df^{n}_{x}(v)\|}{\|v\|} 
      = \lim_{n\to\infty}
\frac{1}{n}\sum_{i=0}^{n-1}\phi(F^{i}(x,v)) 
= \int \phi d\mu.
      \end{equation}
 It remains to show that \( \lambda\in\mathcal L(f) \), i.e. that 
 there actually exists a measure \( \mu^{*}\in\mathcal M(f) \) such 
 that \( \lambda\in\mathcal L(\mu^{*} \)), as opposed to the 
 possibility that the limit in \eqref{lyap2} exists by complete 
 coincidence for some exceptional point.     
 By Lemma \eqref{pistar} the measure \( \mu^{*}=\pi^{*}\mu \) is 
 invariant and ergodic for \( f \). Moreover we have \( 
 \mu^{*}(\pi(A)) = 1 \). Thus, for \( \mu^{*} \) almost every \( x\in 
 \pi(A) \) and every \( v\in T_{x}M \) such that \( (x,v)\in A \) we 
 have that equation \eqref{lyap2} holds and thus \( \lambda \) is one of 
 the Lyapunov exponents associated to the measure \( \mu^{*} \) and 
 therefore belongs to \( \mathcal L(f) \). 
    \end{proof}  
  
\begin{proof}[Proof of Theorem \ref{minimum principle}]
    By Lemma \ref{minimumlemma}, there exists a measure 
    \( \mu\in \mathcal E(F) \) which minimizes the integral \( 
    \int\phi d\mu\). By Lemma \ref{lemmalyap2} there exists a Lyapunov 
    exponent \( \lambda\in\mathcal L(f) \) associated to some measure \( 
    \mu^{*}\in\mathcal E(f) \), which satisfies \( \lambda = \int\phi 
    d\mu \). Thus it just remains to show that there exist no other 
    Lyapunov exponents \( \lambda'<\lambda \).  By Lemma 
    \ref{lemmalyap1} this would imply the existence of a measure \( 
    \mu'\in\mathcal M(F) \) such that \( \int\phi 
    d\mu'=\lambda'<\lambda = \int \phi d\mu \) which contradicts the 
    minimality of \( \int\phi d\mu \) over all measures in \( \mathcal 
    M(F) \). 
 \end{proof}

\section{Uniform expansivity}

We now prove Theorem \ref{uniform}. We first reformulate the 
definition of uniform hyperbolicity in the following clearly 
equivalent form:  
    there exists constants \( \lambda > 0 \) and \( N>0 \) such that 
    \[ 
    \|Df_{x}(v)\|\geq e^{\lambda n} \|v\| \quad \forall \ 
    x\in \Lambda\ v\in T_{x}M, 
    \ n\geq N. 
    \]
\begin{proof}[Proof of Theorem \ref{uniform}]
  By Theorem \ref{minimum principle} we can choose a constant \( 
  \lambda \) satisfying 
    \[ 
    \inf\{\mathcal L(f)\}= \lambda' > \lambda > 0
    \]
    We assume by contradiction that 
    there exists a sequence of times \( n_{k}\to\infty \), 
    a sequence of points \( x_{k} \in \Lambda\) and 
   a sequence of vectors \( v_{k}\in T_{x_{k}}M \)  such that 
   \[
   \|Df^{n_{k}}_{x_{k}}(v_{k})\|< e^{\lambda n_{k}}\|v_{k}\|. 
   \] 
   Now consider the sequence of measures 
   \[ 
   \mu_{n_{k}}= \frac{1}{n_{k}}\sum_{i=0}^{n_{k}-1} 
   \delta_{F^{i}(x_{k},v_{k})}
   \]
   and any (invariant) limit measure \( \mu \) of this sequence. 
   For simplicity 
   we shall suppose without loss of generality that \( 
   \mu_{n_{k}}\to\mu \). Then we have 
   \[ 
   \frac{1}{n_{k}}\sum_{i=0}^{n_{k}-1} \phi (F^{i}(x_{k},v_{k})) = 
   \int \phi d\mu_{n_{k}} < \lambda
   \]
   for every \( n_{k} \). In particular 
   \[ 
    \lim_{k\to\infty}\frac{1}{n_{k}}\sum_{i=0}^{n_{k}-1} 
    \phi (F^{i}(x_{k},v_{k})) = 
   \lim_{k\to\infty}\int \phi d\mu_{n_{k}} =\int\phi d\mu \leq  
   \lambda. 
   \]
However, by Lemma \ref{lemmalyap2} this implies \( \lambda\in\mathcal 
L(f) \) which contradicts our choice of \( \lambda \). 
   \end{proof}

\end{document}